\documentclass[final,12pt]{elsarticle}

\usepackage{amsmath}
\usepackage{amsfonts}
\usepackage{amsthm}

\newtheorem{theorem}{Theorem}[section]

\theoremstyle{remark}
\newtheorem{rem}[theorem]{Remark}

\journal{}

\begin{document}
\begin{frontmatter}

\title{Integral equations and applications}

\author{A. G. Ramm}

\address{ Mathematics Department, Kansas State University,\\
Manhattan, KS 66502, USA\\
email: ramm@math.ksu.edu}


\end{frontmatter}

The goal of  this Section is to formulate some of the basic results
on the theory of integral equations and mention some of its
applications. The literature of this subject is very large. Proofs
are not given due to the space restriction. The results are taken
from the works mentioned in the references.

\section{Fredholm equations} \label{FredholmEquations}
\subsection{Fredholm alternative}
One of the most important results of the theory of integral
equations is the Fredholm alternative. The results in this
Subsection are taken from \cite{KA}, \cite{R419},\cite{R451},
\cite{Kr1}.

Consider the equation
\begin{equation}\label{eq1} u = Ku +
f,\qquad Ku = \int_{D}K(x,y)u(y)dy. \end{equation} Here $f$ and $K(x,y)$ are
given functions, $K(x,y)$ is called the kernel of the operator $K$.
Equation (\ref{eq1}) is considered usually in $H = L^2(D)$ or $C(D)$.
The first is a Hilbert space with the norm $||u|| := \big(\int_D |u|^2
dx\big)^{1/2}$ and the second is a Banach space with the norm $||u|| =
\sup_{x \in D}|u(x)|$. If the domain $D \subset \mathbb{R}^n$ is assumed
bounded, then $K$ is compact in $H$ if, for example,
$\int_D\int_D|K(x,y)|^2dxdy < \infty$, and $K$ is compact in $C(D)$ if,
for example, $K(x,y)$ is a continuous function, or, more generally,
$sup_{x\in D}\int_D |K(x,y)|dy<\infty.$

\noindent For equation (\ref{eq1}) with compact operators, the basic
result is the Fredholm alternative. To formulate it one needs some
notations. Consider equation (\ref{eq1}) with a parameter $\mu \in
\mathbb{C}$:
\begin{equation}\label{eq2} u = \mu Ku + f, \end{equation}
and its adjoint equation
\begin{equation}\label{eq3} v = \bar{\mu}K^{*}v
+ g,\qquad  K^{*}(x,y) = \overline{K(y,x)}, \end{equation}
where the overbar stands
for complex conjugate. By ($2_0$) and ($3_0$) the
corresponding equations with $f = 0$, respectively, $g = 0$ are denoted.
If $K$ is compact so is $K^{*}$. In $H$ one has $(Ku, v) = (u, K^{*}v)$,
where $(u,v) = \int_D u\bar{v}dx$. If equation ($2_0$) has a
non-trivial solution $u_0$, $u_0 \neq 0$, then the corresponding
parameter $\mu$ is called the characteristic value of $K$ and the
non-trivial solution $u_0$ is called an eigenfunction of $K$. The value
$\lambda = \frac{1}{\mu}$ is then called an eigenvalue of $K$. Similar
terminology is used for equation ($3_0$). Let us denote by $N(I -
\mu K) = \{u: u = \mu Ku\}$, $R(I - \mu K) = \{f : f = (I - \mu K)u\}$.
Let us formulate the Fredholm alternative:
\begin{theorem}\label{theorem1} If $N(I - \mu K) = \{0\}$, then $N(I -
\bar{\mu}K^{*}) = \{0\}, R(I - \mu K) = H, R(I - \bar{\mu}K^{*}) =
H$.$\\$ If $N(I - \mu K) \neq \{0\}$, then $\dim N(I - \mu K) = n <
\infty, \dim N(I - \bar{\mu} K^{*}) = n, f \in R(I - \mu K)$ if and only
if  $f \perp N(I - \bar{\mu}K^{*})$, and $g \in R(I -
\bar{\mu}K^{*})$ if and only if $g \perp N(I - \mu K)$. \end{theorem}

\begin{rem}\label{remark1} Conditions $f \perp N(I - \bar{\mu}K^{*})$
and $g \perp N(I - \mu K)$ are often written as $(f, v_{0j}) = 0, 1 \leq
j \leq n$, and, respectively, $(g, u_{0j}) = 0, 1 \leq j \leq n$, where
$\displaystyle\{v_{0j}\}_{j = 1}^n$ is a basis of the subspace $N(I -
\bar{\mu}K^{*})$, and $\displaystyle\{u_{0j}\}_{j = 1}^n$ is a basis of
the subspace $N(I - {\mu}K)$. \end{rem}

\begin{rem}\label{remark2} Theorem \ref{theorem1} is very useful in
applications. It says that equation (\ref{eq1}) is solvable for any $f$
as long as equation ($2_0$) has only the trivial solution $v_0 =
0$. A simple proof of Theorem \ref{theorem1} can be found in
\cite{R419} and \cite{R451}. \end{rem}

\begin{rem}\label{remark3} Theorem \ref{theorem1} can be stated
similarly for equations (\ref{eq2}) and (\ref{eq3}) in $C(D).$ \end{rem}

\begin{rem}\label{remark4} If $n > 0$ and $f \in R(I - \mu K)$, then the
general solution to equation (\ref{eq2}) is $u = \tilde{u} +
\displaystyle\sum_{j = 1}^n c_ju_{0j}$, where $c_j, 1 \leq j \leq n$,
are arbitrary constants and $\tilde{u}$ is a particular solution to
(\ref{eq2}). Similar result holds for $v.$ \end{rem}

\subsection{Degenerate kernels} Suppose that \begin{equation}\label{eq4}
K(x,y) = \sum_{m = 1}^Ma_m(x)\overline{b_m(y)}, \end{equation} where
the functions $a_m$ and $b_m$ are linearly independent, $a_m, b_m \in H$.
The operator $K$ with the degenerate kernel (\ref{eq4}) is called a
finite-rank operator. Equation (\ref{eq2}) with degenerate kernel
(\ref{eq4}) can be reduced to a linear algebraic system. Denote
\begin{equation}\label{eq5} u_m = \int_D u\overline{b_m}dx.
\end{equation} Multiply equation (\ref{eq2}) by $\overline{b_p}$,
integrate over $D$, and get \begin{equation}\label{eq6} u_p = \mu\sum_{m
= 1}^Ma_{pm}u_m + f_p,\quad  1 \leq p \leq M, \end{equation} where $a_{pm} :=
(a_m, b_p)$. There is a one-to-one correspondence between solutions to
equation (\ref{eq2}) and (\ref{eq6}): \begin{equation}\label{eq7} u(x) =
\mu\sum_{m = 1}^Ma_m(x)u_m + f(x). \end{equation}
{\it Exercise
1.} Solve the equation $u(x) = \mu\displaystyle\int_0^{\pi}\sin(x -
t)u(t)dt + 1$.

{\it Hint:} Use the formula $\sin(x-t) = \sin x \cos t - \cos x \sin
t$. Find characteristic values and eigenfunctions of the operator with
the kernel $\sin(x - t)$ on the interval $(0, \pi)$.

\begin{rem}\label{remark5} If $K$ is a compact integral operator in
$L^2(D)$ then there is a sequence $K_n(x, y)$ of degenerate kernels such
that $\displaystyle\lim_{n \to \infty}||K - K_n|| = 0$, where $||K||$ is
the norm of the operator $K$, $||K|| := \displaystyle\sup_{u \neq
0}\frac{||Ku||}{||u||}$. \end{rem}

\subsection{Volterra equations} If $K(x,y) = 0$ for $y > x$ and $\dim D
= 1$ then equation (\ref{eq2}) is called Volterra equation. It is of the
form \begin{equation}\label{eq8} u(x) = \mu\int_a^xK(x,y)u(y)dy + f(x),
\quad a \leq x \leq b. \end{equation} This is a particular case of the
Fredholm equations. But the Volterra equation (\ref{eq8}) has special
property. \begin{theorem}\label{theorem2} Equation (\ref{eq8}) is
uniquely solvable for any $\mu$ and any $f$ in $C(D)$, if the function
$K(x,y)$ is continuous in the region $a \leq y \leq x,\quad a\leq x \leq b$.
The solution to (\ref{eq8}) can be obtained by iterations:
\begin{equation}\label{eq9} u_{n+1}(x) = \mu Vu_n + f, \quad u_1 = f;\quad
 Vu :=
\int_a^xK(x,y)u(y)dy. \end{equation} \end{theorem}

\subsection{Selfadjoint operators} If $K^* = K$ then $K$ is called
selfadjoint. In this case equation (\ref{eq2}) has at least one
characteristic value, all characteristic values are real numbers, and the
corresponding eigenfunctions are orthogonal. One can solve equation
(\ref{eq2}) with selfadjoint operator $K$ by the formula
\begin{equation}\label{eq10} u = \sum_{j = 1}^\infty \frac{\mu_j
f_j}{\mu_j - \lambda}u_{0j}(x), f_j = (f, u_{0j}), \end{equation} where
$u_{0j} = \mu_jKu_{oj}, \,\, j = 1,2, \dots$.$\\$ For the eigenvalues of a
selfadjoint compact operator $K$ the following result (minimax
representation) holds. We write $K \geq 0$ if $(Ku,u) \geq 0, \forall u
\in H$, and $A \leq B$ if $(Au,u) \leq (Bu,u), \forall u \in H.$
\begin{theorem}\label{theorem3} If $K = K^* \geq 0$ is compact in $H$,
$Ku_j = \lambda_ju_j,\,\,\, u_j \neq 0$ then
\begin{equation}\label{eq11} \lambda_{j + 1} = \max_{u \perp
L_j}\frac{(Ku,u)}{||u||^2},\quad j = 0,1,\dots \end{equation} Here
$L_0 = \{0\}, L_j$ is the subspace with the basis $\{u_1, \dots,
u_j\},\,\, ||u_j|| = 1,\,\, (u_i, u_m) = \delta_{im}$. Also
\begin{equation}\label{eq12} \lambda_{j + 1} = \min_M\max_{u \perp
M}\frac{(Ku,u)}{||u||^2}, \quad \dim M = j, \end{equation} where $M$
runs through the set of $j$-dimensional subspaces of $H$.
\end{theorem} One has $\lambda_1 \geq \lambda_2 \geq \dots \geq 0$. The
set of all eigenfunctions of $K$, including the ones corresponding
to the zero eigenvalue, if such an eigenvalue exists, forms an
orthonormal basis of $H$.
\begin{rem}\label{remark6} If $0 \leq A \leq B$ are compact operators in
$H$, then $\lambda_j(A) \leq \lambda_j(B).$ \end{rem}
\begin{rem}\label{remark7} If $A$ is a linear compact operator in $H$,
then $A = U|A|$, where $|A| = (A^*A)^{1/2} \geq 0$ and $U$ is a partial
isometry, $$U: R(|A|) \to R(A),\,\, ||U|A|f|| = ||Af||\,\,\, \forall f \in
H,\,\,
 N(U)= N(A),$$
and $\overline{R(|A|)} = \overline{R(A^*)},$ where the overline
denotes
the closure in $H$. One has \begin{equation}\label{eq13} Au = \sum_{j =
1}^{r(A)}s_j(A)(u,u_j)v_j. \end{equation} Here $r(A) = \dim R(|A|)$, $s_j
\geq 0$ are the eigenvalues of $|A|$. They are called the $s$-values of
$A$, $\{u_j\}$ and $\{v_j\}$ are orthonormal systems, $|A|u_j = s_ju_j,$
$s_1 \geq s_2 \geq \dots$, $v_j = Uu_j.$ \end{rem}
\begin{rem}\label{remark8} If $A$ is a linear compact operator in $H$
then \begin{equation}\label{eq14} s_{j + 1}(A) = \min_{K \in K_j}||A -
K||, \end{equation} where $K_j$ is the set of finite rank operators of
dimension $\leq j$, that is, $dim R(K_j)=j$. \end{rem}

\subsection{Equations of the first kind}
If $Ku = \lambda u + f$ and $\lambda = 0$, then one has
\begin{equation}\label{eq15} Ku = f. \end{equation} If $K$ is compact in
$H$, then equation (\ref{eq15}) is called Fredholm equation of the first
kind. The linear set $R(K)$ is not closed unless $K$ is a finite-rank
operator. Therefore, small pertubation of $f \in R(K)$ may lead to an
equation which has no solutions, or, if it has a solution, this solution
may differ very much from the solution to equation (\ref{eq15}). Such
problems are called ill-posed. They are important in many applications,
especially in inverse problems, \cite{R470}, \cite{I} . The usual statement of
the problem of solving an ill-posed equation (\ref{eq15}) consists of
the following. It is assumed that equation (\ref{eq15}) is solvable for
the exact data $f$,
possibly non-uniquely, that the exact data $f$ are not known but the
``noisy'' data $f_\delta$ are known, and $||f_\delta - f|| \leq \delta$,
where $\delta>0$ is known.
Given $f_\delta$, one wants to find $u_\delta$ such that
$\displaystyle\lim_{\delta \to 0}||u - u_\delta|| = 0.$

There are many methods to solve the above problem: the DSM (Dynamical
Systems Method) \cite{R499} \cite{R612}, iterative methods,
variational regularization, and their variations, see
\cite{I} and \cite{R470}.

Examples of ill-posed problems of practical interest include stable
numerical differentiation, stable summation of the Fourier series and
integrals with perturbed coefficients, stable solution to linear
algebraic systems with large condition numbers, solving Fredholm and
Volterra integral equations of the first kind, various deconvolution
problems, various inverse problems, the Cauchy problem for the Laplace
equation and other elliptic equations, the backward heat problem,
tomography, and many other problems, see \cite{R470}.

\subsection{Numerical solution of Fredholm integral equations}
Let us
describe the projection method for solving equation (\ref{eq2}) with
compact operator $K$. Let $L_n \subset H$ be a linear subspace of $H$,
$\dim L_n = n$. Assume throughout that the sequence $L_n \subset L_{n +
1}$ is limit dense in $H$, that is, $\displaystyle\lim_{a \to
\infty}\rho(f, L_n) = 0$ for any $f \in H$, where $\rho(f, L_n) =
\displaystyle\inf_{u \in L_n}||f - u||$. Let $P_n$ be an orthogonal
projection operator on $L_n$, that is, $P_n^2 = P_n,$
$P_nf \in L_n$, $ ||f -
P_nf|| \leq ||f - u||, \forall u \,\, in \,\,L_n.$

Assume that $N(I - K) = \{0\}$. Then $||(I - K)^{-1}|| \leq c$. Consider
the projection method for solving equation (\ref{eq1}):
\begin{equation}\label{eq16} P_n(I - K)P_nu_n = P_nf. \end{equation}
\begin{theorem}\label{theorem4} If $N(I - K) = \{0\}$, $K$ is compact,
and $\{L_n\}$ is limit dense, then equation (\ref{eq16}) is uniquely
solvable for all $n$ sufficiently large, and $\displaystyle\lim_{n \to
\infty}||u_n - u|| = 0$, where $u$ is the solution to equation
(\ref{eq1}). \end{theorem} \begin{rem}\label{remark9} Projection method
is valid also for solving equations $Au = f$, where $A: X \to Y$ is a
linear operator, $X$ and $Y$ are Hilbert spaces. Let $Y_n \subset Y$ and
$X_n \subset X$ be $n$-dimensional subspaces, $P_n$ and $Q_n$ are
projection operators on $X_n$ and $Y_n$ respectively. Assume that the
sequences $\{X_n\}$ and $\{Y_n\}$ are limit dense and consider the
projection equation \begin{equation}\label{eq17} Q_nAP_nu_n = Q_nf,\quad u_n
\in X_n. \end{equation} \end{rem} \begin{theorem}\label{theorem5} Assume
that $A$ is a bounded linear operator, $||A^{-1}|| \leq c$, and
\begin{equation}\label{eq18}
||Q_nAP_nu|| \geq c||P_nu||,\quad  \forall u \in X,
\end{equation} where $c > 0$ stands for various positive constants
independent of $n$ and $u$. Then equation (\ref{eq17}) is uniquely
solvable for any $f \in Y$ for all sufficiently large $n$, and
\begin{equation}\label{eq19} \lim_{n \to \infty}||u - u_n|| = 0.
\end{equation} Conversely, if equation (\ref{eq17}) is uniquely solvable
for any $f \in Y$ for all sufficiently large $n$, and (\ref{eq19})
holds, then (\ref{eq18}) holds. \end{theorem}
\begin{rem}\label{remark10} If (\ref{eq18}) holds for $A$, then it will
hold for $A + B$ if $||B|| > 0$ is sufficiently small. If (\ref{eq18})
holds for $A$, then it holds for $A+B$ if $B$ is compact and $||(A +
B)^{-1}|| \leq c.$ \end{rem} Iterative methods for solving integral
equation (\ref{eq1}) are popular. If $||K|| < 1$, then equation
(\ref{eq1}) has a unique solution for every $f \in H$, and this solution
can be obtained by the iterative method \begin{equation}\label{eq20}
u_{n + 1} = Ku_n + f,\quad  \lim_{n \to \infty}||u_n - u|| = 0, \end{equation}
where $u_1 \in H$ is arbitrary.

Define the spectral radius of a linear bounded operator $A$ in a Banach
space by the formula $\rho(A) = \displaystyle\lim_{n \to
\infty}||A^n||^{1/n}$. It is known that this limit exists. If $a$ is a
number, then $\rho(aA) = |a|\rho(A)$. The basic fact is:

{\it If $|\lambda| >
\rho(A)$, then the equation $\lambda u = Au + f$ is uniquely solvable
for any $f$, and the iterative process $\lambda u_{n + 1} = Au_n + f$
converges to $u = (\lambda I - A)^{-1}f$ for any $f \in H$ and any
initial approximation $u_1$. }

Given an equation $Au = f$, can one
transform it into an equation $u = Bu + g$ with $\rho(B) < 1,$ so that
it can be solved by iterations? To do so, one may look for an operator
$T$ such that $(A + T)u = Tu + f, (A + T)^{-1}$ is bounded, so $u = (A +
T)^{-1}Tu + (A + T)^{-1}f$, and $\rho((A + T)^{-1}T) < 1$. Finding such
an operator $T$ is, in general, not simple. If
$A = A^*,\,\,\, 0 \leq m \leq A \leq M$, where
$m$ and $M$ are constants, one may transform equation $Au = f$ to the
equivalent form $u = u - \frac{2Au}{m + M} + \frac{2}{m + M}f$, and get
$\rho(I - \frac{2A}{m + M}) \leq \frac{M - m}{M + m} < 1.$ If $A \neq
A^*$ then equation $A^*Au = A^*f$ is useful since $A^*A$ is a
selfadjoint operator. The above results can be found, for example, in
\cite{Kr1} and \cite{Z}.

Equation $Au = f$ is often convenient to consider in a cone of a
Banach space $X$, see \cite{Z}. A cone is a closed convex set $K
\subset X$ which contains with a point $u$ all the points $tu, t
\geq 0$, and such that $u \in K$ and $-u \in K$ imply $u = 0$. One
writes $u \leq v$ if $v - u \in K$. This defines a semiorder in $X$.
If $K$ is invariant with respect to $K$, that is, $u \in K$ implies
$Au \in K$, then one can give conditions for the solvability of the
equation $Au = f$ in $K$. For example, let $X = C(D), D \subset
\mathbb{R}^m$ is a bounded domain, $K = \{u: u \geq 0\},$
$$Au =
\int_DA(x,y)u(y)dy, \quad A(x,y) \geq 0.$$

The following result is useful:

 {\it Suppose that there exist functions
$v$ and $w$ in $K$ such that $Av \geq v$, $Aw \leq w$, and $ v \leq
w$. Let us assume that every monotone bounded sequence in $K$, $u_1
\leq u_2 \leq \dots \leq U$ has a limit: $\exists u \in K$ such that
$\displaystyle\lim_{n \to \infty}||u_n - u|| = 0$. Then there exists
$u \in K$ such that $u = Au$, where $ v \leq u \leq w.$}

 A cone is
called normal if $0 \leq u \leq v$ implies $||u|| \leq c||v||,$
where $ c> 0$ is constant independent of $u$ and $v$.

\section{Integral equations with special kernels} \subsection{Equations
with displacement kernel}

The results in this Subsection are taken from \cite{G}, \cite{Z}.

Let us mention some equations with displacement kernels. We start
with the equation
$$\lambda u(t) - \int_0^tK(t -
s)u(s)ds = f(t),\quad t \geq 0.$$ This equation is solved by taking Laplace
transform: $\lambda \bar{u} - \bar{K}\bar{u} = \bar{f}$, $\bar{u} =
(\lambda - \bar{K})^{-1}\bar{f},$ where  $\bar{u} := Lu :=
\int_0^{\infty}e^{-pt}u(t)dt,$ and
$$ u(t) = L^{-1}\bar{u} = \frac{1}{2\pi
i}\displaystyle\int_{\sigma - i\infty}^{\sigma +
i\infty}e^{pt}\bar{u}(p)dp.$$
Equation
$$\lambda u(t) -
\int_{-\infty}^{\infty}K(t - s)u(s)ds = f(t),\quad -\infty < t < \infty,$$ is
solved by taking the Fourier transform: $\tilde{u}(\xi) :=
\int_{-\infty}^{\infty}u(t)e^{-i\xi t}dt.$ One has $\tilde{u}(\xi) = (\lambda -\tilde{K}(\xi))^{-1}\tilde{f},$ and $$ u(t) =
\frac{1}{2\pi}\int_{-\infty}^{\infty}\tilde{u}(\xi)e^{i\xi t}d\xi.$$ In
these formal calculations one assumes that $1 - \tilde{K}(\xi) \neq 0$,
see \cite{G}, \cite{Mu}.

More complicated is Wiener-Hopf equation \begin{equation}\label{eq20a}
u(t) - \int_0^\infty K(t - s)u(s)ds = f(t),\quad  t \geq 0. \end{equation}
Assume that $K(t) \in L(-\infty, \infty)$ and $1 - \tilde{K}(\xi) \neq
0,\,\, \forall \xi \in (-\infty, \infty)$. Define the index: $\kappa :=
-\displaystyle\frac{1}{2\pi}\arg[1 -
\tilde{K}(\xi)]|_{-\infty}^{\infty}$. The simplest result is:

{\it If and only if $\kappa = 0$ equation (\ref{eq20a}) is uniquely solvable in
$L^p(-\infty, \infty)$.}

The solution can be obtained analytically by the following scheme.
Denote $u_+(t) := \begin{cases} u(t), &t \geq 0 \\ 0, & t <
0\end{cases}.$ Write equation (\ref{eq20a}) as $$u_+(t) =
\int_{-\infty}^\infty K(t - s)u_+(s)ds + u_-(t) + f_+,\quad -\infty
< t < \infty,$$ where $u_-(t) = 0$ for $t \geq 0$ is an unknown
function. Take the Fourier transform of this equation to get $[1 -
\tilde{K}(\xi)]\tilde{u}_+ = \tilde{f}_+ + \tilde{u}_-$. Write $1 -
\tilde{K}(\xi) = \tilde{K}_+(\xi)\tilde{K}_-(\xi)$, where
$\tilde{K}_+ $ and $\tilde{K}_-$ are analytic in the upper,
respectively, lower half-plane functions. Then
$$\tilde{K}_+(\xi)\tilde{u}_+(\xi) =
\tilde{f}_+(\xi)\tilde{K}_-^{-1}(\xi) +
\tilde{u}_-(\xi)\tilde{K}_-(\xi).$$
Let $P_+$ be the projection operator
on the space of functions  analytic in the upper half-plane.  Then
$$\tilde{u}_+ =
\tilde{K}_+^{-1}(\xi)P_+(\tilde{f}_+(\xi)\tilde{K}_-^{-1}(\xi)).$$
If $\tilde{u}_+$ is known, then $u_+(t)$ is known, and the Wiener-Hopf
equation is solved. The non-trivial part of this solution is
the factorization of the function $1-\tilde{K}(\xi)$.

\subsection{Equations basic in random fields estimation theory}
In this Subsection the results from \cite{R486} are given.

Let $u(x) = s(x) + n(x), x \in \mathbb{R}^2, s(x)$ and $n(x)$ are
random fields, $s(x)$ is a ``useful signal'' and $n(x)$ is
``noise''. Assume that the mean values $\overline{s(x)} =
\overline{n(x)} = 0$, the overbar denotes mean value. Assume that
covariance functions
$$\overline{u^*(x)u(y)} := R(x,y),\quad  \overline{u^*(x)s(y)} := f(x,y),$$ are
known. The star here denotes complex conjugate. Suppose that a
linear estimate $Lu$ of the observations in a domain $D$ is $$Lu =
\int_D h(x,y)u(y)dy,$$ where $h(x,y)$ is, in general, a
distributional kernel. Then the optimal estimate of $As$, where $A$
is a linear operator, is found by solving the optimization problem
\begin{equation}\label{eq21} \epsilon := \overline{(Lu - As)^2} =
\min . \end{equation} A necessary condition for the kernel (filter)
$h(x,y)$ to solve this problem is the integral equation:
\begin{equation}\label{eq22} Rh := \int_D R(x,y)h(z,y)dy = f(x,z),
\quad x,z \in \bar{D} := D \cup \Gamma,
\end{equation} where $\Gamma := \partial D$ is the boundary of $D$.
Since $z$ enters into equation (\ref{eq22}) as a parameter, one has
to study the following integral equation
\begin{equation}\label{eq23} \int_D R(x,y)h(y)dy = f(x), \quad x \in
\bar{D}. \end{equation} The questions of interest are:
\begin{itemize}
 \item[1)] Under what assumptions equation (\ref{eq23}) has a
unique solution in some space of distributions?
 \item[2)] What is the order of singularity of this solution?
 \item[3)] Is this solution stable under small (in some sense)
perturbations of $f$? Does it solve estimation problem (\ref{eq21})?
 \item[4)] How does one calculate this solution analytically and
numerically? \end{itemize}
If a distribution $h = D^lh_1$, where $h_1
\in L^2_{\mbox{loc}}$, then the order of singularity of $h$,
$\mbox{ordsing}h$, is $|l|$ if $h_1$ is not smoother than
$L^2_{\mbox{loc}}$. By $|l|$ one denotes $l_1 + l_2 + \dots + l_r, D^l =
\displaystyle\frac{\partial^{|l|}}{\partial x_1^{l_1} \dots \partial
x_r^{l_r}}$. For the distribution theory one can consult, for example,
\cite{GS}.

Let us define a class of kernels $R(x,y)$, or, which is the same, a
class $\mathcal{R}$ of integral equations (\ref{eq23}), for which the
questions 1) - 4) can be answered.

Let $L$ be a selfadjoint elliptic operator of order $s$ in
$L^2(\mathbb{R}^r)$, $\Lambda$, $d\rho(\lambda)$ and
$\Phi(x,y,\lambda)$ are, respectively, its spectrum, spectral
measure and spectral kernel. This means that the spectral function
$E_\lambda$ of the operator $L$ has the kernel $E_\lambda(x,y) =
\int_\Lambda \Phi(x,y,\lambda)d\rho(\lambda)$ and a function
$\phi(L)$ has the kernel $$\phi(L)(x,y) = \int_\Lambda
\phi(\lambda)\Phi(x,y,\lambda)d\rho(\lambda),$$ where
$\phi(\lambda)$ is an arbitrary function such that
$\int_\Lambda|\phi(\lambda)|^2d(E_\lambda f,f) < \infty$ for any $f
\in D(\phi(L))$.

Define class $\mathcal{R}$ of the kernels $R(x,y)$ by the formula
\begin{equation}\label{eq24} R(x,y) = \int_\Lambda
\frac{P(\lambda)}{Q(\lambda)}\Phi(x,y,\lambda) d\rho(\lambda),
\end{equation} where $P(\lambda) > 0$ and $Q(\lambda) > 0$ are
polynomials of degree $p$ and, respectively, $q$, and  $p \leq q$.
Let $\alpha := 0.5(q - p)s$. The number $\alpha \geq 0$ is an
integer because $p$ and $q$ are even integers if the polynomials
 $P> 0$ and $Q> 0$ are positive for all values of $\lambda$.
 Denote by $H^{\alpha} = H^{\alpha}(D)$ the Sobolev space and
by $\dot{H}^{-\alpha}$ its dual space with respect to the $L^2(D)$
inner product. Our basic results are formulated in the following
theorems, obtained in  \cite{R486}.
\begin{theorem}\label{theorem6} If $R(x,y) \in \mathcal{R}$, then the
operator $R$ in equation (\ref{eq23}) is an isomorphism between the
spaces $\dot{H}^{-\alpha}$ and $H^{\alpha}$. The solution to
equation (\ref{eq23}) of minimal order of singularity,
$\mbox{ordsing} h \leq \alpha$, exists, is unique, and can be
calculated by the formula
\begin{equation}\label{eq25} h(x) = Q(L)G, G(x) = \begin{cases} g(x) +
v(x) &\text{in } D \\ u(x) &\text{in } D' := \mathbb{R}^r \setminus D
\end{cases}
\end{equation}
where $g(x) \in H^{0.5s(p + q)}$ is an arbitrary fixed solution to
the equation \begin{equation}\label{eq26} P(L)g = f \text{ in }
D,\quad f\in H^\alpha, \end{equation} and the functions $u(x)$ and
$v(x)$ are the unique solution to the problem
\begin{equation}\label{eq27} Q(L)u = 0 \text{ in } D',\quad  u(\infty) = 0,
\end{equation} \begin{equation}\label{eq28} P(L)v = 0 \text{ in } D,
\end{equation} \begin{equation}\label{eq29} \partial_N^ju =
\partial_N^j(v + g) \text{ on } \partial D,\,\,\, 0 \leq j \leq 0.5s(p +
q) -
1. \end{equation} Here $\partial_N u$ is the derivative of $u$ along
the outer unit normal $N$.
\end{theorem}

Theorem \ref{theorem6} gives answers to questions 1) - 4), except
the question of numerical solution of equation (\ref{eq23}). This
question is discussed in \cite{R486}.
\begin{rem}\label{remark11} If one considers the operator $R \in
\mathcal{R}$ as an operator in $L^2(D)$,  $R: L^2(D) \to L^2(D)$ and
denotes by $\lambda_j = \lambda_j(D)$ its eigenvalues, $\lambda_1
\geq \lambda_2 \geq \dots \geq 0$, then \begin{equation}\label{eq30}
\lambda_j = cj^{-(q - p)s/r}[1 + o(1)] \quad \text{ as }\,\, j \to
\infty,
\end{equation} where $c = \gamma^{(q - p)s/r},$  $ \gamma :=
(2\pi)^{-r}\int_D \eta(x)dx$, and the $\eta(x)$ is defined in
\cite{R486}. \end{rem} {\it Example.} Consider equation (\ref{eq23})
with $D = [-1,1]$, $R(x,y) = \exp(-|x - y|)$, $L = -i\frac{d}{dx}$,
$r = 1$, $s = 1$, $P(\lambda) = 1$, $Q(\lambda) = (\lambda^2 +
1)/2$,  $ \Phi(x,y,\lambda) = (2\pi)^{-1}\exp \{i\lambda(x - y)\}$,
and $d\rho(\lambda = d\lambda)$. Then
$$\int_{-1}^1\exp(-|x -y|)h(y)dy=f,\quad -1 \leq x \leq 1,$$
$\alpha = 1$,  and Theorem \ref{theorem6} yields the following
formula for the solution $h$:
$$h(x) = \frac{-f'' + f}{2} + \frac{f'(1) + f(1)}{2}\delta(x - 1) +
\frac{-f'(-1) + f(-1)}{2}\delta(x + 1),$$ where $\delta$ is the
delta-function.

\section{Singular integral equations} \subsection{One-dimensional
singular integral equations}
The results in this Section are taken from \cite{MP}, \cite{Sh} and
\cite{Z}.

Consider the equation
\begin{equation}\label{eq31} Au := a(t)u(t) +
\frac{1}{i\pi}\int_L\frac{M(t,s)}{s - t}u(s)ds = f(t).
\end{equation} Here the functions $a(t),\, M(t,s),$ and $ f(t)$
are known.
Assume that they satisfy the H\H{o}lder condition. The contour $L$
is a closed smooth curve on the complex plane, $L$ is the boundary
of a connected smooth domain $D$, and $D'$ is its complement on the
complex plane. Since $M(t,s)$ is assumed H\H{o}lder-continuous, one
can transform equation (\ref{eq31}) to the form
\begin{equation}\label{eq32} Au := a(t)u +
\frac{b(t)}{i\pi}\int_L\frac{u(s)ds}{s - t} + \int_L K(t,s)u(s)ds =
f(t), \end{equation} where the operator with kernel $K(t,s)$ is of
Fredholm type, $b(t) = M(t,t)$, $K(t,s) = \frac{1}{2\pi}\frac{M(t,s)
- M(t,t)}{s - t}$. The operator $A$ in \eqref{eq32} is of the form
$A = A_0 + K$, where $A_0u = a(t)u + Su$, and $Su =
\frac{1}{i\pi}\int_L\frac{u(s)ds}{s - t}$. The singular operator $S$
is defined as the limit $\displaystyle\lim_{\epsilon \to
0}\frac{1}{i\pi}\int_{|s - t| > \epsilon}\frac{u(s)ds}{s - t}$.

The operator $S$ maps $L^p(L)$ into itself and is bounded, $1 < p <
\infty$. Let us denote the space of H\H{o}lder-continuous functions
on $L$ with exponent $\gamma$ by $\mbox{Lip}_\gamma(L)$. The
operator $S$ maps the space $\mbox{Lip}_\gamma(L)$ into itself and
is bounded in this space. If $\Phi(z) = \frac{1}{2\pi
i}\int_L\frac{u(s)ds}{s - z}$ and $ u \in \mbox{Lip}_\alpha(L)$,
then
\begin{equation}\label{eq33} \Phi^{\pm}(t) = \pm \frac{u(t)}{2} +
\frac{1}{2\pi i}\int_L\frac{u(s)ds}{s - t},
\end{equation} where $+ (-)$ denotes the limit $z \to t,\, z \in D \,(D').$

If $u \in L^p(L), 1 < p < \infty$, then \begin{equation}\label{eq34}
\frac{d}{dt}\int_Lu(s)\ln\frac{1}{|s - t|}ds = i\pi u(t) +
\int_L\frac{u(s)ds}{s - t}. \end{equation} One has $S^2 = I$, that
is, for $1 < p < \infty$, \begin{equation}\label{eq35}
\frac{1}{(i\pi)^2}\int_L\frac{1}{\tau - t} \bigg(
\int_L\frac{u(s)ds}{s - \tau}ds \bigg)d\tau = u(t),\quad  \forall u
\in L^p(L). \end{equation} Let $n = \dim N(A)$ and $n^* := \dim
N(A^*)$. Then index of the operator $A$ is defined by the formula
$\mbox{ind}A := n - n^*$ if at least one of the numbers $n$ or $n^*$
is finite.

If $A=B+K$, where $B$ is an isomorphism and $K$ is compact, then $A$
is Fredholm operator and $\mbox{ind} A = 0$. If $A$ is a singular
integral operator \eqref{eq32}, then its index may be not zero,
\begin{equation} \label{eq36}
    \kappa:=\mbox{ind} A=\frac{1}{2\pi} \int_L d \arg
\frac{a(t)-b(t)}{a(t)+b(t)}, \end{equation} where $ \arg f$ denotes the
argument of the function $f$. The operator $A$ with non-zero index
is called a Noether operator. Formula \eqref{eq36} was derived
by F.Noether in 1921. Equation \eqref{eq32} with $K(x,y)=0$ can be
reduced to the Riemann problem. Namely, if $\Phi(z)=\frac{1}{2\pi
i}\int_L \frac{u(s)}{s-z}ds$, then $u(t)=\Phi^+(t)-\Phi^-(t)$ and
$\Phi^+(t)+\Phi^-(t)=Su$ by formula \eqref{eq33}.

Thus, equation \eqref{eq32} with $K(x,y)=0$ can be written as
\begin{equation*}
    a(t)(\Phi^+-\Phi^-)+b(t)(\Phi^+ +\Phi^-)=f,
\end{equation*}
or as
\begin{equation} \label{eq37}
    \Phi^+(t)=\frac{a(t)-b(t)}{a(t)+b(t)}\Phi^- + \frac{f}{a(t)+b(t)}, \quad t \in L.
\end{equation}
This problem is written usually as
\begin{equation} \label{eq38}
    \Phi^+(t)=G(t)\Phi^-(t)+g(t), \quad t \in L,
\end{equation}
and is a Riemann problem for finding piece-wise analytic function
$\Phi(z)$ from the boundary condition \eqref{eq38}.

 Assume that
$a(t)\pm b(t)\neq 0, \forall t \in L$, and that $L$ is a closed
smooth curve, the boundary of a simply connected bounded domain $D$.

If one solves problem \eqref{eq38} then the solution to equation
 \eqref{eq32} (with $K(x,y)=0$)
is given by the formula $u=\Phi^+ - \Phi^-$. Let $\kappa=\mbox{ind}
G$.
\begin{theorem} \label{thm7}
    If $\kappa \geq 0$ then problem \eqref{eq38} is solvable for
    any $g$ and its general solution is
    \begin{equation} \label{eq39}
        \Phi(z)=\frac{X(z)}{2\pi i}\int_L \frac{g(s)ds}{X^+(s)(s-z)}+X(z)P_\kappa(z),
    \end{equation}
    where
    \begin{align}
        &X^+(z)=exp(\Gamma_+(z)),\quad X^-(z)=z^{-\kappa}exp(\Gamma_-(z)), \label{eq40}\\
        &\Gamma(z)=\frac{1}{2\pi i} \int_L \frac{\ln[s^{-\kappa}G(s)]}{s-z}ds, \label{eq41}
    \end{align}
    where $P_\kappa(z)$ is a polynomial of degree $\kappa$ with arbitrary
    coefficients. If $\kappa=-1$, then problem \eqref{eq38} is
    uniquely solvable and $X^-(\infty)=0$. If $\kappa <-1$
    then problem \eqref{eq38}, in general, does not have a solution.
    For its solvability, it is necessary and sufficient that
    \begin{equation} \label{eq42}
        \int_L \frac{g(s)}{X^+(s)}s^{k-1}ds=0, k=1,2,\ldots, -\kappa-1.
    \end{equation}
    If these conditions hold, then the solution to problem \eqref{eq38}
    is given by formula \eqref{eq39} with $P_\kappa(z)\equiv 0$.
\end{theorem}
This theorem is proved, for example, in \cite{G} and \cite{MP}. It
gives also an analytic formula for solving integral equation
\eqref{eq32} with $k(x,y)=0$.

If one consider equation \eqref{eq32} with $k(x,y)\neq0$, then this
equation can be transformed to a Fredholm-type equation by inverting
the operator $A_0$. A detailed theory can be found in \cite{G} and
\cite{Mu}.

\subsection{Multidimensional singular integral equations}
Consider the singular integral
\begin{equation} \label{eq43}
    Au:=\int_{\mathbf{R}^m} r^{-m}f(x,\theta)u(y)dy, \quad x,y \in \mathbf{R}^m,\quad r=|x-y|,
\end{equation}
Assume that $u(x)=O(\frac{1}{|x|^k})$ for large $|x|$, that $k>0,$ and
$f$ is a bounded function continuous, with respect to $\theta$. Let
 $\int_S f(x,\theta)dS=0$, where $S$ stands in this section for the
unit sphere in $\mathbf{R}^m$.

Suppose that $|u(x)-u(y)|\leq c_1 r^\alpha (1+|x|^2)^{-k/2}$ for $ r
\leq 1$, where $ c_1=const$,  and $|u(x)|\leq c_2(1+|x|^2)^{-k/2}$,
where $c_2=const$. Denote this set of functions $u$ by
$A_{\alpha,k}$. If the first inequality holds, but the second is
replaced by the inequality $|u(x)| \leq c_3
(1+|x|^2)^{-k/2}\ln(1+|x|^2)$, then this set of $u$ is denoted
$A'_{\alpha,k}$.
\begin{theorem} \label{thm8}
    The operator \eqref{eq43} maps $A_{\alpha,k}$ with $k \leq m$
    into $A'_{\alpha,k}$, and $A_{\alpha,k}$ with $k > m$ into $A_{\alpha,m}$.
\end{theorem}
Let $D \subset \mathbf{R}^m$ be a domain, possibly $D =
\mathbf{R}^m$, the function $f(x,\theta)$ be continuously
differentiable with respect to both variables, and $u$ satisfies the
Dini condition, that is $\sup_{||x-y||<\delta;x,y\in D}
|u(x)-u(y)|=w(u,t)$, where $\int_0^\delta \frac{w(u,t)}{t}dt <
\infty, \quad \delta=const>0$.

Then
\begin{equation} \label{eq44}
    \frac{\partial}{\partial x_j} \int_D \frac{f(x,\theta)}{r^{m-1}}u(y)dy=
    \int_D u(y) \frac{\partial}{\partial x_j}\left[\frac{f(x,\theta)}
    {r^{m-1}}\right]dy-u(x)\int_S f(x,\theta)\cos(r_{x_j})ds.
\end{equation}
Assume that
\begin{equation} \label{eq45}
    \int_S |f(x,\theta)|^{p'} dS \leq C_0=const,
    \quad \frac{1}{p}+\frac{1}{p'}=1,
\end{equation}
and let
\begin{equation} \label{eq46}
    A_0 u =\int_{\mathbf{R}^m} \frac{f(x,\theta)}{r^{m}} u(y)dy.
\end{equation}
\begin{theorem} \label{thm9}
    If \eqref{eq45} holds then
    $$||A_0||_{L^p} \leq c\sup_x ||f(x,\theta)||_{L^{p'}(S)}.$$
\end{theorem}
The operator \eqref{eq46} is a particular case of a
pseudodifferential operator defined by the formula
\begin{equation} \label{eq47}
    A u = \int_{\mathbf{R}^m} \int_{\mathbf{R}^m}
    e^{i(x-y)\cdot \xi} \sigma_A(x,\xi)u(y)dyd\xi,
\end{equation}
where the function $\sigma_A(x,\xi)$ is called the symbol of $A$.
For the operator $Au=a(x)u(x)+A_0 u$, where $A_0$ is defined in
formula \eqref{eq46}, the symbol of $A$ is defined as
$\sigma(x,\xi)=a(x)+\widetilde{K}(x,\xi)$, where
$K(x,x-y)=\frac{f(x,\theta)}{r^m}$, and $\widetilde{K}(x,\xi)$ is
the Fourier transform of $K(x,z)$ with respect to $z$. One has
$$\sigma_{A+B}(x,\xi)=\sigma_A(x,\xi)+\sigma_B(x,\xi),\quad
\sigma_{AB}(x,\xi)=\sigma_A(x,\xi)\sigma_B(x,\xi).$$

It is possible to estimate the norm of the operator \eqref{eq46}
in terms of its symbol,
\begin{equation*}
    \sigma_{A_0}(x,\xi)=\int_{\mathbf{R}^m} \frac{f(x,y^0)}
    {|y|^m} e^{{-iy \cdot z}} dy, \quad y^0=\frac{y}{|y|}, \quad \xi=\frac{z}{|z|}.
\end{equation*}
If $T$ is a compact operator in $L^2(\mathbf{R}^m)$ then its symbol
is equal to zero. Let $Au=a(x)u+A_0u$, and
$\sigma_A(x,\xi)=a(x)+\sigma_{A_0}(x,\xi)$. If $\sigma_A(x,\xi)$ is
sufficiently smooth and does not vanish, then the singular integral
equation $Au+Tu=f$ is of Fredholm type and its index is zero. In
general, index of a system of singular integral equations (and of
more general systems) is calculated in the work \cite{AZ}.
Additional material about singular integral equations and
pseudodifferential operators one finds in \cite{MP} and \cite{Sh}.

\section{Nonlinear integral equations}
The results from this Subsection are taken from \cite{Kr}, \cite{Z},
\cite{KZ}, \cite{KA}.

Let us call the operators $$Uu=\int_D K(x,t)u(t)dt, \quad Hu=\int_D
K(x,t)f(t,u(t))dt,$$
 and $Fu:=f(t,u(t))$, respectively, Urysohn, Hammerstein, and Nemytskii operators.

The operator $F$ acts from $L^p:=L^p(D)$ into $L^q:=L^q(D)$ if and
only if  $$|f(t,u)| \leq a(t)+b|u|^{p/q},\quad p<\infty,\quad a(t)
\in L^q,\quad b=const.$$ If $p=\infty$, then $F: L^\infty\to L^q$ if
and only if $|f(t,u)| \leq a_h(t)$, $a_h(t) \in L^q$, and $
|u(t)|\leq h, \,\,0 \leq h < \infty$. The function $f(t,u),\, t \in
D,\,\, u \in (-\infty, \infty)$, is assumed to satisfy the
Caratheodory conditions, that is, this function is continuous with
respect to $u$ for almost all $t \in D$ and is measurable with
respect to $t \in D$ for all $u \in (-\infty, \infty)$. If $F: L^p
\to L^q$ and $f(t,u)$ satisfies the Caratheodory conditions, then
$F$ is compact provided that one of the following conditions holds:
1) $q<\infty$,\,\,  2) $q=\infty,\,\, p < \infty,\,\, f(t,u)\equiv
a(t)$,\,\,  3) $q=p=\infty,\,\, |f(t,u)-f(t,v)|\leq \phi_h(u-v)$
with respect to $z$ for any $h>0$.

If $f$ satisfies the Caratheodory conditions, and $K(x,t)$ is a
measurable function on $D \times D,$ and $F: L^p \to L^r,\,\, K: L^r
\to L^q$, then the operator $H=KF: L^p \to L^q$, $H$ is continuous
if $r<\infty$, and $H$ is compact if $r<\infty$ and $K$ is compact.
These results one can find, for example, in \cite{Z}.

A nonlinear operator $A:X \to Y$ is Fr\'echet differentiable if
$$A(u+h)-A(u)=Bh+w(u,h),$$ where $B: X \to Y$ is a linear bounded
operator and $\lim_{||h||_X \to 0} \frac{||w(u,h)||_Y}{||h||_X}=0$.

Consider the equation
\begin{equation} \label{eq48}
    u(x)=\mu \int_D K(x,t,u(t))dt+f(x):=Au,
\end{equation}
where $D \subset \mathbf{R}^m$ is a bounded closed set, meas$D>0$,
the functions $K(x,t,u)$ and $f(x)$ are given, the function $u$ is
unknown, and $\mu$ is a number.
\begin{theorem} \label{thm10}
    Let us assume that $A$ is a contraction on a set $M$, that is,
$||A(u)-A(v)\le q||u-v||$, where
     $u,v \in M$, $M$ is a subset of a Banach space $X$, $0<q<1$ is a number,
      and $A: M \to M$. Then equation \eqref{eq48} has a unique solution
      $u$ in $M$, $u=\lim_{n \to \infty}u_n$ , where $u_{n+1}=A(u_n), u_0 \in M$,
      and $||u_n-u|| \leq \frac{q^n}{1-q}||u_1-u_0||, n=1,2,\ldots$.
\end{theorem}
This result is known as the {\it contraction mapping principle}.
\begin{theorem} \label{thm11}
    Assume that $A$ is a compact operator, $M \subset X$ is a
    bounded closed convex set, and $A: M \to M$. Then equation
    \eqref{eq48} has a solution.
\end{theorem}
This result is called {\it the Schauder principle}. The solution is a
fixed point of the mapping $A$, that is, $u=A(u)$. Uniqueness of the
solution in Theorem \ref{thm11} is not claimed: in general, there
can be more than one solution.
\begin{rem}
    A version of Theorem \ref{thm11} can be formulated as follows:
\end{rem}
\begin{theorem} \label{thm13}
    Assume that $A$ is a continuous operator which maps a convex closed
    set $M$ into its compact subset. Then $A$ has a fixed point in $M$.
\end{theorem}
Let us formulate {\it the Leray-Schauder} principle:
\begin{theorem} \label{thm14}
    Let $A$ be a compact operator, $A: X \to X$, $X$ be a Banach
    space. Let all the solutions $u(\lambda)$ of the
    equation $u=A(u,\lambda),\,\, 0 \leq \lambda \leq 1$,
    satisfy the estimate $\sup_{0 \leq \lambda \leq 1} ||u(\lambda)|| \leq a < \infty$.
    Let $||A(u,0)|| \leq b$ for $b>a$ and $||u||=b$. Then the equation
    $u=A(u,1)$ has a fixed point in the ball $||u||\leq a$.
\end{theorem}
Let $D \subset X$ be a bounded convex domain, $A^nu \subset
\overline{D}$ if $u\in \overline{D}, n=1,2,\ldots$, $\overline{D}$
is the closure of $D$, $A^nu \neq u$ for $n>n_0$ and $u \in \partial
D$, where $\partial D$ is the boundary of $D$. Finally, let us
assume that $A$ is compact.
\begin{theorem} \label{thm15}
    Under the above assumptions $A$ has a fixed point in $D$.
\end{theorem}
This theorem is proved in \cite{KZ}

Consider the equation $u=H(u), H: L^2(D) \to L^2(D)$, and assume
that $K(x,t)$ is positive-definite, continuous kernel, $f(t,u)$
is continuous with respect to $t \in D$ and $u \in (-\infty, \infty)$
function, such that $\int_0^u f(t,s) ds \leq 0.5 au^2+b$,
where $a < \Lambda$, and $\Lambda$ is the  maximal eigenvalue of
the linear operator $K$ in $L^2(D)$. Then there exists a fixed point of the operator $H$.
This result can be found in \cite{Z}.

\section{Applications}
\subsection{General remarks.}
The results in this Subsection are taken from \cite{R190},
\cite{R595} and \cite{R632}.

There are many applications of integral equations. In this section,
there is no space to describe in detail applications to solving
boundary value problems by integral equations involving potentials
of single and double layers, applications of integral equations
in the elasticity theory (\cite{Mu}, \cite{Ku}),
applications to acoustics electrodynamics, etc.

We will restrict this Subsection to some questions not considered by
 other authors. The first question deals with the possibility
 to express potentials of the single layer by potentials of the
 double layer and vice versa. The results are taken from \cite{R190}.
   The second application deals
 with the wave scattering by small bodies of an arbitrary shape.
 These results are taken from \cite{R595} and \cite{R632}.

\subsection{Potentials of a single and double layers.}
\begin{equation} \label{eq49}
    V(x)=\int_S g(x,t) \sigma(t) dt, \quad g(x,t)=\frac{e^{ik|x-t|}}{4\pi |x-t|},
\end{equation}
where $S$ is a boundary of a smooth bounded domain $D \subset
\mathbf{R}^3$ with the boundary $S$, $\sigma(t) \in Lip_\gamma
(S),\,\, 0 < \gamma \leq 1$.

It is known that $V(x) \in C(\mathbf{R}^3), V(\infty)=0$,
\begin{equation} \label{eq50}
    V_N^\pm(s)=\frac{A\sigma\pm\sigma(s)}{2}, \quad A\sigma:=2\int_S \frac{\partial g(s,t)}
    {\partial N_s}\sigma(t)dt,
\end{equation}
where $N=N_s$ is the unit normal to $S$ at the point $s$ pointing
out of $D$, $+(-)$ denotes the limiting value of the normal
derivative of $V$ when $ x \to s \in S$ and $x \in D\,\, (D')$,
where $D':=\mathbf{R}^3 \setminus D$.

Potential of double layer is defined as follows:
$$W(x):=\int_S \frac{\partial g(x,t)}{\partial N_t}\mu(t)dt.$$
One has the following properties of $W$:
\begin{align}
    &W^\pm(s)=W(s)\mp\frac{\mu(t)}{2}, \quad W(s):=\int_S \frac{\partial g(s,t)}{\partial N_t}\mu(t)dt \label{eq51} \\
    &W^+_{N_s}=W^-_{N_s}. \label{eq52}
\end{align}
\begin{theorem} \label{thm16}
    For any $V(W)$ there exists a $W(V)$ such that $W=V$ in $D$. The $V(W)$ is uniquely defined.
\end{theorem}
\begin{theorem} \label{thm17}
    A necessary and sufficient condition for $V(\sigma)=W(\mu)$ in $D'$ is
    \begin{equation} \label{eq53}
        \int_S Vh_j dt = 0, \quad 1 \leq j \leq r', \quad (I+A)h_j=0,
    \end{equation}
    where the set $\{h_j\}$ forms a basis of $N(I+A)$, and $A$ is defined in \eqref{eq50}.

    A necessary and sufficient condition for $W(\mu)=V(\sigma)$ in $D'$ is:
    \begin{equation} \label{eq54}
        \int_S W\sigma_j dt = 0, \quad 1 \leq j \leq r, \quad A\sigma_j-\sigma_j=0,
    \end{equation}
    where the set $\{\sigma_j\}$ forms a basis of $N(I-A)$.
\end{theorem}
Theorem \ref{thm16} and \ref{thm17} are proved in \cite{R190}.

\subsection{Wave scattering by small bodies of an arbitrary shape.} \label{sec5.3}
 The results in this Subsection are taken from \cite{R595} and \cite{R632}.
 Consider the wave scattering problem:
\begin{align}
    &(\nabla^2+k^2)u=0 \text{ in } D':=\mathbf{R}^3\setminus D, \quad u_N=\zeta u \text{ on }S=\partial D,\label{eq55} \\
    &u=u_0+v, \quad u_0=e^{ik\alpha\cdot x}, \quad \alpha \in S^2, \label{eq56} \\
    &\frac{\partial v}{\partial r}-ikv=o\left(\frac{1}{r}\right),\quad r=|x|\to\infty. \label{eq57}
\end{align}
Here $k=const>0$, $S^2$ is the unit sphere in $\mathbf{R}^3, \zeta=const$
is a given parameter, the boundary impedance, $D$ is a small
body, a particle,  $S$ is its boundary, which we assume H\"older-continuous,
and $N$ is the unit normal to $S$ pointing out of $D$.
If $Im\zeta \leq 0$ then problem \eqref{eq55}-\eqref{eq57}
has exactly one solution, $v$ is the scattered field,
\begin{equation} \label{eq58}
    v=A(\beta, \alpha, \kappa)\frac{e^{ikr}}{r}+o\left(\frac{1}{r}\right),\quad |x|=r \to \infty, \quad \beta=\frac{x}{r},
\end{equation}
$A(\beta,\alpha, \kappa)$ is called the scattering amplitude.

{\it Our basic assumption is the smallness of the body $D$:
this body is small if $ka\ll 1$, where $a=0.5diam D$.}

In applications $a$
is called the characteristic size of $D$.

Let us look for the solution to problem \eqref{eq55}-\eqref{eq57} of the form
\begin{equation} \label{eq59}
    u(x)=u_0(x)+\int_S g(x,t) \sigma(t)dt, \quad
g(x,t)=\frac{e^{ik|x-t|}}{4\pi|x-t|},
\end{equation}
and write
\begin{equation} \label{eq60}
    \int_S g(x,t) \sigma(t)dt=g(x,x_1)Q+\int_S[g(x,t)-g(x,x_1)]\sigma(t)dt,
\end{equation}
where $x_1 \in D$ is an arbitrary fixed point inside $D$, and
\begin{equation} \label{eq61}
    Q:=\int_S \sigma(t)dt.
\end{equation}
If $ka\ll 1$ and $|x-x_1|:=d\gg a$, then
\begin{equation} \label{eq62}
    \left|\int_S[g(x,t)-g(x,x_1)]\sigma(t)dt\right| \ll |g(x,x_1)|Q.
\end{equation}
Therefore the scattering problem \eqref{eq55}-\eqref{eq57} has an
approximate solution of the form:
\begin{equation} \label{eq63}
    u(x)=u_0(x)+g(x,x_1)Q, \quad ka\ll 1, \quad |x-x_1|\gg a.
\end{equation}
The solution is reduced to finding just one number $Q$
in contrast with the usual methods, based on the boundary
integral equation for the unknown function $\sigma(t)$.

To find the main term of the asymptotic of $Q$ as $a \to 0$, one
uses the exact boundary integral equation:
\begin{equation} \label{eq64}
    u_{0N}(s)-\zeta u_0(s)+\frac{A\sigma-\sigma}{2} -
\zeta\int_S g(s,t)\sigma(t)dt=0, \quad s\in S,
\end{equation}
where formula \eqref{eq50} was used.

We do not want to solve equation \eqref{eq64} which is only
numerically possible, but want to derive an asymptotically exact
analytic formula for $Q$ as $a \to 0$. Integrate both sides of
formula \eqref{eq64} over $S$. The first term is equal to
 $\int_S u_{0N}ds = \int_D \nabla^2 u_0 dx \simeq O(a^3)$.
 The sign $\simeq$ stands for the equality up to the terms of higher order
 of smallness as $a\to 0$. The second term is equal to
 $-\zeta\int_S u_0 ds \simeq -\zeta u_0(x_1)|S|$, where
 $|S|=O(a^2)$ is the surface area of $S$. The third term is
 equal to $\frac{1}{2}\int_S A\sigma ds-\frac{1}{2}Q$. When $a \to 0$
 one checks that $A\sigma \simeq A_0 \sigma$, where $A_0=A|_{k=0}$,
 and $\int_S A_0\sigma dt=-\int_S \sigma dt$. Therefore the third term
 is equal to $-Q$ up to the term of higher order of smallness.
 The fourth term is equal to $-\zeta\int_S dt\sigma(t) \int_S g(s,t)ds=o(Q)$,
 as $a \to 0$. Thus,
\begin{equation} \label{eq65}
    Q\simeq -\zeta|S|u_0(x_1), \quad a \to 0,
\end{equation}
and $u_0(x_1)\simeq 1$,
\begin{equation} \label{eq66}
    A(\beta,\alpha,\kappa)=-\frac{\zeta |S| u_0(x_1)}{4\pi}.
\end{equation}
For the scattering by a single small body one can choose $x_1\in D$
to be the origin, and take $u_0(x_1)=1$. But in the many-body
scattering problem the role of $u_0$ is played by the effective field
which depends on $x$.
Formulas \eqref{eq63}-\eqref{65} solve the scattering problem \eqref{eq55}-\eqref{eq57} for one
 small body $D$ of an arbitrary shape if the impedance  boundary condition
 \eqref{eq55} is imposed. The scattering in this case is isotropic and
 $A(\beta,\alpha,\kappa)=O(a^2\zeta)$.

If the boundary condition is the Dirichlet one, $u|_S =0$, then one derives that
\begin{equation} \label{eq67}
    Q \simeq -C u_0, \quad a\to 0; \quad u_0 \simeq 1,
\end{equation}
and
\begin{equation} \label{eq68}
    A(\beta,\alpha,\kappa)=-\frac{C}{4\pi}, \,\,\,a \to 0,\quad u_0 \simeq 1, \quad C=O(a).
\end{equation}
Here $C$ is the electrical capacitance of the perfect conductor with the shape $D$.

In this case the scattering is isotropic and $|A|=O(a), a \to 0$. Thus, the scattered
field is much larger than for the impedance boundary condition.

For the Neumann boundary condition $u_N=0$ on $S$ the scattering is
anisotropic and $A(\beta,\alpha,\kappa)=O(a^3)$, which is much
smaller than for impedance boundary condition.

One has for the Neumann boundary condition the following formula:
\begin{equation} \label{eq69}
    A(\beta,\alpha,\kappa)=\frac{D}{4\pi}\left(ik\beta_{pq}
    \frac{\partial u_0}{\partial x_q}\beta_p+\nabla^2 u_0\right),
\end{equation}
where $\beta_{pq}$ is some tensor, $|D|$ is the volume of $|D|$,
 $\beta_p=\lim_{|x|\to\infty}\frac{x_p}{x}$, $x_p:=x\cdot e_p$ is the
 $p-$th coordinate of a vector $x\in \mathbf{R}^3$,  $u_0$ and its derivatives
 are calculated at an arbitrary point inside $D$. This point can be
 chosen as the origin of the coordinate system. Since $D$ is small,
 the choice of this point does not influence the results.
The tensor $\beta_{pq}$ is defined for a body with volume $V$ and
boundary $S$ as follows: $$\beta_{pq}=\frac{1}{V}\int_S t_p \sigma_q
dt,$$ where $\sigma_q$ solves the equation
$\sigma_q=A\sigma_q-2N_q$, $N_q:=N\cdot e_q$, and $\{e_q\}_{q=1}^3$
is a Cartesian basis of $\mathbf{R}^3$.


Finally, consider the many-body scattering problem. Its statement can be
written also as \eqref{eq55}-\eqref{eq57} but now $D=\cup_{m=1}^M D_m$
is the union of many small bodies, $\zeta=\zeta_m$ on $S_m=\partial D_m$.

One looks for the solution of the form
\begin{equation} \label{eq70}
    u(x)=u_0(x) +\sum_{m=1}^M \int_{S_m} g(x,t) \sigma_m(t) dt.
\end{equation}
Let us define the effective field in the medium by the formula
\begin{equation} \label{eq71}
    u_e(x)=u_0(x) +\sum_{m\neq j} g(x,x_m)Q_m, \quad |x-x_m|\sim a.
\end{equation}
Here $x_m \in D_m$ are arbitrary fixed points, and $|x-x_m|\sim a$
means that $|x-x_m|$ is of the size of $ a$.

Assume that these points are distributed by the formula:
\begin{equation} \label{eq72}
    \mathcal{N}(\Delta)=\frac{1}{a^{2-\kappa}} \int_\Delta N(x)dx[1+o(1)], \quad a\to 0.
\end{equation}
Here $\mathcal{N}(\Delta)=\sum_{x_m \in \Delta} 1$ is the number of
points in an arbitrary open set $\Delta$, $N(x)\geq 0$ is a given
function, $0\leq \kappa<1$ is a parameter, and an experimentalist
can choose $N(x)$ and $\kappa$ as he wishes. Let us assume that
$\zeta_m = \frac{h(x_m)}{a^\kappa}$, where $h(x)$ is an arbitrary
continuous in $D$ function such that $Im h(x) \leq 0$. This function
can be chosen by an experimenter as he wishes.

Under these assumptions in \cite{R632} the following results are proved.

Assume that $|S_m|=ca^2$, where $c>0$ is a constant
depending on the shape of $D_m$. Then $Q_m \simeq -c h(x_m)u_e(x_m)a^{2-\kappa}$.
 Denote $u_e(x_m):=u_m, h(x_m):=h_m$. Then one can find the unknown numbers
 $u_m$ from the following linear algebraic system (LAS):
\begin{equation} \label{eq73}
    u_j=u_{0j} -c\sum_{m\neq j}^M g_{jm} h_m u_m a^{2-\kappa}, \quad 1 \leq j \leq M,\quad
    g_{jm}:=\frac{e^{ik|x_j-x_m|}}{4\pi|x_j-x_m|}.
\end{equation}
LAS \eqref{eq73} can be reduced to a LAS of much smaller order. Namely, denote by $\Omega$
a finite domain in which all the small bodies are located. Partition $\Omega$
into a union of $P$ cubes $\Delta_p$ with the side $b=b(a)$, the cubes are
non-intersecting in the sense that they do not have common interior points,
 they can have only pieces of common boundary. Denote by $d=d(a)$ the
 smallest distance between neighboring bodies $D_m$. Assume that
\begin{equation} \label{eq74}
    a \ll d \ll b, \quad ka\ll 1.
\end{equation}
Then system \eqref{eq73} can be reduced to the following LAS:
\begin{equation} \label{eq75}
    u_q=u_{0q} -c\sum_{p\neq q}^P g_{qp} h_p u_p N_p |\Delta_p|, \quad 1 \leq q \leq P,
\end{equation}
where $N_p:=N(x_p)$, $x_p \in \Delta_p$,  is an arbitrary point,
$N(x)$ is the function from formula \eqref{eq72}, and $|\Delta_p|$
is the volume of the cube $\Delta_p$.

The effective field has a limit as $a\to 0$, and this limit solves the equation
\begin{equation} \label{eq76}
    u(x)=u_{0}(x) -c\int_\Omega g(x,y) h(y) N(y)u(y)dy,
\end{equation}
where the constant $c$ is the constant in the definition $|S_m|=ca^2$.
One may consider by the same method the case of small bodies of various
sizes. In this case $c=c_m$, but we do not go into details.
If the small bodies are spheres of radius $a$ then $c=4\pi$.

Applying the operator
 $\nabla^2+k^2$ to the equation \eqref{eq76} one gets
\begin{equation} \label{eq77}
    (\nabla^2+k^2n^2(x))u:=(\nabla^2+k^2-c N(y) h(y))u=0 \text{ in } \mathbf{R}^3.
\end{equation}
Therefore, embedding of many small impedance particles, distributed according
formula \eqref{eq72} leads to a medium with a zero refraction coefficient.
\begin{equation} \label{eq78}
    n^2(x)=1-k^{-2}c N(x)h(x).
\end{equation}
Since $N(x)$ and $h(x)$ can be chosen as one wishes, with the only restrictions
$N(x)\geq 0, Im h(x)\leq 0$, one can create a medium with a desired
refraction coefficient by formula \eqref{eq78} if one chooses suitable
$N(x)$ and $h(x)$.

A similar theory is developed for electromagnetic wave scattering by small
impedance bodies in \cite{R628}.


\end{document}